\documentclass[12pt]{article}

\usepackage[colorlinks=true, pdfstartview=FitV, linkcolor=blue, citecolor=blue, urlcolor=blue]{hyperref}

\usepackage{amssymb,amsmath,amscd,array}
\usepackage{times, verbatim}
\usepackage{graphicx}
\usepackage{graphics}
\usepackage{epic}

\DeclareFontFamily{OT1}{rsfs}{}
\DeclareFontShape{OT1}{rsfs}{n}{it}{<-> rsfs10}{}
\DeclareMathAlphabet{\mathscr}{OT1}{rsfs}{n}{it}
\let\mathcal=\mathscr

\usepackage{amssymb,amsmath, amscd}
\usepackage{times, verbatim}
\usepackage{graphicx}

\DeclareFontFamily{OT1}{rsfs}{}
\DeclareFontShape{OT1}{rsfs}{n}{it}{<-> rsfs10}{}
\DeclareMathAlphabet{\mathscr}{OT1}{rsfs}{n}{it}

\DeclareMathOperator{\Gal}{Gal}

\DeclareMathOperator{\SO}{SO}

\DeclareMathOperator{\SL}{SL}
\DeclareMathOperator{\PGL}{PGL}
\DeclareMathOperator{\GL}{GL}
\DeclareMathOperator{\Tr}{Tr}

\DeclareMathOperator{\Ind}{Ind}

\DeclareMathOperator{\ad}{ad}
\DeclareMathOperator{\art}{art}
\DeclareMathOperator{\sw}{sw}
\DeclareMathOperator{\Fr}{Fr}
\DeclareMathOperator{\ZZ}{Z}

\def\R{{\mathbb R}}

\def\bC{{\mathbb C}}
\def\Q{{\mathbb Q}}
\def\Z{{\mathbb Z}}
\def\A{{\mathbb A}}

\title{Simple supercuspidals and the Langlands correspondence}

\author{Benedict H Gross}

\begin{document}
\maketitle

\tableofcontents

\section{Formal degrees}

In the fall of $2007$, Mark Reeder and I were thinking about the formal degrees of discrete series representations of $p$-adic groups $G$. The formal degree is a generalization of the dimension of an irreducible representation, when $G$ is compact. It depends on the choice of a Haar measure $dg$ on $G$: if the irreducible representation $\pi$ is induced from a finite dimensional representation $W$ of an open, compact subgroup $K$, then the formal degree of $\pi$ is given by $\deg(\pi) = \dim(W)/\int_K dg$.

Mark had determined the formal degree in several interesting cases \cite{R1}, including the series of depth zero supercuspidal representations that he had constructed with Stephen DeBacker \cite{DR}. We came across a paper by Kaoru Hiraga, Atsushi Ichino, and Tamotsu Ikeda, which gave a beautiful conjectural formula for the formal degree in terms of the $L$-function and $\epsilon$-factor of the adjoint representation of the Langlands parameter  \cite[\S1]{HII}. This worked perfectly in the depth zero case \cite[\S3.5]{HII}, and we wanted to test it on more ramified parameters, where the adjoint $L$-function is trivial.

Naturally, we started with the group $\SL_2(\Q_p)$. In the course of our work, Mark found a new construction of irreducible representations, which we called simple supercuspidals. His discovery, when viewed through the prism of the Langlands correspondence, has led to some interesting mathematics, which I want to survey here.

\section{The conjecture}

The formula that Hiraga, Ichino, and Ikeda propose for the formal degree of a discrete series representation $\pi$ depends on the (still conjectural) Langlands parametrization of irreducible representations. We will state this when ${\bf G}$ is a split, simply-connected group defined over the ring of integers $A$ of a $p$-adic field $k$, and $G = {\bf G}(k)$ is the locally compact group of points, with open compact subgroup $K = {\bf G}(A)$. Let $F$ be the residue field of $A$, which is finite of order $q$. As in \cite{Gr2} and \cite[\S1]{HII}, we choose the Haar measure $dg$ on $G$ which gives the open compact subgroup $K$ volume 
$ \int_K dg = \#{\bf G}(F) / q^{\dim G}.$

The Langlands parameter $\phi$ of an irreducible complex representation of $G$ is a homomorphism
$$\phi: W(k) \times \SL_2(\bC) \rightarrow \hat{G}(\bC)$$
where $W(k)$ is the Weil group of $k$ and $\hat{G}$ is the dual group. This satisfies a number of conditions, which are detailed in \cite[\S 3.2]{GR}. In particular, we consider two homomorphisms equivalent if they are conjugate in $\hat{G}(\bC)$. For simplicity, we will only define all the terms in their conjecture for the formal degree when the Langlands parameter is trivial on $\SL_2(\bC)$. 

The composition of $\phi$ with the adjoint representation $\hat{\frak g}$ of the dual group gives a complex, orthogonal representation of the Weil group $W(k) \rightarrow \SO(\hat{\frak g})$. We recall that $W(k) = I.\Fr^{\Z}$, where $I$ is the inertia subgroup of the Galois group of $k$ and $\Fr$ is a geometric Frobenius. The $L$-function of the adjoint representation is defined by the formula
$$L(\phi, \ad, s) = \det(1 - \Fr.q^{-s}|\hat{\frak g}^I)^{-1}.$$
Associated to a parameter $\phi$ is the centralizer $C_{\phi}$ of the image in $\hat{G}$; for parameters of representations in the discrete series, this centralizer is conjectured to be finite. 
Under this assumption, the Weil group $W(k)$ has no invariants on $\hat{\frak g}$ and the adjoint $L$-function is regular at $s=0$.

 We let $\art(\ad)$ be the Artin conductor of the adjoint representation of $W(k)$. This is a non-negative integer, of the form \cite[\S 2.2]{GR}
$$\art(\ad) = \dim(\hat{\frak g}/\hat{ \frak g}^I) + \sw(\ad).$$
Here $\sw(\ad)$ is the Swan conductor, which is non-zero if and only if the parameter is wildly ramified. The $\epsilon$-factor of the adjoint representation has the form \cite[\S 2.2]{GR}
$$\epsilon(\phi, \ad,s) = \pm q^{\art(\ad) (1/2 - s)}$$
The sign $ \pm = \epsilon(\phi,\ad, 1/2)$ will not concern us here.

Every discrete Langlands parameter should correspond to a finite $L$-packet of irreducible representations $\pi$ of $G$, and elements of the $L$-packet should be in bijection with the set of irreducible representations $\rho$ of the centralizer $C_{\phi}$.Their conjecture for the formal degree a the discrete series representation $\pi$ with Langlands parameter $(\phi, \rho)$, with respect to the Haar measure $dg$ specified above, is
$$\deg(\pi) = \pm (\dim \rho)/\#C_{\phi} \times L(\phi, \ad, 1)\epsilon(\phi, \ad, 0) / L(\phi, \ad, 0). $$

\section{Depth zero for $\SL_2$}

Depth zero supercuspidal representations are all induced from maximal compact subgroups of $G$. They have Langlands parameters which are tamely ramified, with $\dim(\hat{\frak g})^I = \dim (\hat{\frak t})$ equal to the rank $r(G)$. This makes it relatively easy to check that the formal degree conjecture is correct \cite[\S 3.5]{HII}.We will do this exercise (in a special case) for the group $\SL_2$ over $k = \Q_p$, where the Langlands parametrization is known.

When $G = \SL(V) \cong \SL_2(\Q_p)$, where $V$ is a two dimensional vector space over $\Q_p$,  every maximal compact subgroup is the stabilizer of a $\Z_p$-lattice $L$ in $V$. Let $W$ be an irreducible complex representation which is in the discrete series of $\SL(L/pL) \cong \SL_2(\Z/p\Z)$. This finite group is a quotient of the compact group $K = \SL(L) \cong \SL_2(\Z_p)$ and one defines the representation $\pi = \Ind(K|G, W)$ of $G$ by compact induction. The resulting infinite dimensional representation is irreducible and supercuspidal -- its matrix coefficients have compact support, so $\pi$ appears as a submodule of the functions with compact support on $G$. 

There are two orbits of $G$ on the set of lattices in $V$, represented by $L_1$ and $L_2$. Here $L_1$ is any lattice in $V$ and $L_2$ is any sublattice of index $p$ in $L_1$. The stabilizers $K_1 = \SL(L_1)$ and $K_2 = \SL(L_2)$ represent the two conjugacy classes of maximal compact subgroups of $G$, and each has quotient isomorphic to $\SL_2(\Z/p\Z)$. When $p$ is odd, there are two half discrete series representations $W_1$ and $W_2$ of the finite group $\SL_2(\Z/p\Z)$ which have dimension $(p-1)/2$ and are exchanged by an outer automorphism \cite[\S 2]{Gr}. This gives us four irreducible representations $\pi_{ij} = \Ind(K_i|G, W_j)$ of $G$, which all have the same formal degree. The Haar measure $dg$ defined above gives each maximal compact subgroup volume $\#\SL_2(\Z/p\Z)/p^3 = (1 - p^{-2})$. Hence the formal degree, defined using this Haar measure, is given by 
$$\deg(\pi_{ij}) = (p-1)/2(1-p^{-2}) = p^2/2(p+1)$$

The four representations $\pi_{ij}$ all have the same Langlands parameter $\phi: W(\Q_p) \rightarrow \PGL_2(\bC)$, so lie in the same $L$-packet. Their parameter factors through the abelianized Weil group $W(\Q_p)^{ab} \cong \Q_p^*$, and even through its quotient $\Q_p^*/\Q_p^{*2}$ which is a Klein $4$-group. The homomorphism $\phi$ injects this quotient into the dual group $\PGL_2(\bC)$. The image is self-centralizing, and is normalized by a subgroup isomorphic to $S_4$, which permutes the three elements of order $2$ triply transitively. Hence the Langlands parameter $\phi$ is well-defined up to conjugacy, and its centralizer $C_{\phi}$ has order $4$. It is impossible to say which representation $\pi_{ij}$ corresponds to which representation $\rho$ of dimension $1$ of the group $C_{\phi}$ without choosing some additional data (like a generic character of the unipotent radical of $\SL_2$) but they all give the ratio $(\dim \rho)/\#C_{\phi} = 1/4$ in the formal degree conjecture.

The parameter $\phi$ is tamely ramified. The image of the inertia group has order $2$ and fixes a one dimensional subspace of the adjoint representation $\frak{pgl}_2$. Hence the Artin conductor of the adjoint representation is equal to $2$. Frobenius acts through a non-trivial quadratic character on the one dimensional subspace fixed by inertia, so the $L$-function and $\epsilon$-factor of the adjoint representation are given by
$$L(\phi, \ad, s) = (1 + p^{-s})^{-1}~~~~~~\epsilon(\phi, \ad, s) = \pm p^{1 - 2s}.$$
The prediction for the formal degree of each representation in the $L$-packet is 
$$\pm1/4 \times L(\phi, \ad, 1)\epsilon(\phi, \ad, 0) / L(\phi, \ad, 0) = 1/4 \times (1 + p^{-1})^{-1}p/ 2^{-1} = p^2/2(p+1)$$
which checks!

\section{Simple supercuspidals for $\SL_2$} 

When the Langlands parameter $\phi: W(\Q_p) \rightarrow \PGL_2(\bC)$ for a representation of $\SL_2(\Q_p)$ is wildly ramified, it is possible for the inertial invariants  on the adjoint representation $\frak{pgl}_2$ to be trivial. Then the $L$-function $L(\phi, \ad, s) = 1$ and the Artin conductor of the adjoint representation is greater than or equal to $4$, with equality holding when the Swan conductor is equal to $1$. When the Artin conductor is $4$, we have $\epsilon(\phi, \ad, s) = \pm p^{2-4s}$ and the prediction for the formal degree involves the ratio
$$L(\phi, \ad, 1)\epsilon(\phi, \ad, 0) / L(\phi, \ad, 0) = \pm p^2.$$
One possible construction of a discrete series with such a simple formal degree would be to induce a one dimensional representation $\chi$ from a pro-$p$-Sylow subgroup $P$ of $K = \SL(L)$. The subgroup $P$ has index $p^2 -1$ in $K$, and hence has volume $1/p^2$ for our fixed Haar measure. It reduces modulo $p$ to a maximal unipotent subgroup $U$ in $\SL_2$ over $\Z/p\Z$.

The only problem was that it seemed extremely unlikely that the induction of a one dimensional representation of $P$ would produce anything interesting. For example, the only one dimensional representation of $K$ is the trivial representation, and inducing this to $G$ gives a module of infinite length, containing all of the $K$-unramified representations! The same situation occurs for the induction of the trivial character of $P$, or any character of $P$ which factors through the finite quotient $U(\Z/p\Z)$. Were there any other characters of order $p$ that we could use?

I asked Mark if he knew the Frattini quotient of $P$ -- the maximal quotient which is an elementary abelian $p$ group. He worked this out when $p > 2$ (which we will assume in this section) using the theory of affine roots. The pro-$p$-group $P$ is a normal subgroup, of index $p-1$, in the Iwahori subgroup $I = \SL(L_1) \cap \SL(L_2)$. The Frattini quotient of $P$ is isomorphic to $(\Z/p\Z)^2$, and the quotient $I/P = (\Z/p\Z)^*$ acts on this quotient via the two simple affine roots $\{\pm\alpha\}$. Hence the action of $I$ gives a natural decomposition of the Frattini quotient into the sum of two lines. Mark reported, to his enormous surprise, that if the character $\chi$ of $P$ was non-trivial when restricted to each eigenspace, the compactly induced representation $\Ind(P|G, \chi)$ was the direct sum of two irreducible representations $\pi_1$ and $\pi_2$ of $G = \SL_2(\Q_p)$, which were distinguished by their central characters! He showed me a proof using Mackey's irreducibility criterion for an induced representation \cite[\S 7.4]{Se}, using the fact that the double cosets of $I$ in $G$ are indexed by elements in the affine Weyl group. The formal degrees of $\pi_1$ and $\pi_2$ are both $p^2/2$. 

To check that these degrees are consistent with the conjecture of Hiraga, Ichino, and Ikeda, we need to identify the Langlands parameters of $\pi_1$ and $\pi_2$, which are homomorphisms $W(\Q_p) \rightarrow \PGL_2(\bC)$. Unlike the depth zero case described above, this homomorphism does not factor through the abelianized Weil group -- the image in $\PGL_2(\bC)$ is a dihedral group of order $2p$. We'll describe these parameters in section $6$.

\section{Simple supercuspidals in the split case}

Mark found these irreducible representations for $\SL_2(\Q_p)$ in early February of $2008$, and soon saw how to generalize his construction to split, simply-connected, almost simple groups $G$ over a $p$-adic field $k$. (For the general case, see \cite[\S2.6]{RY}.) Let $A$ be the ring of integers of $k$ and let $F$ be the residue field, of order $q$. Let $I$ be an Iwahori subgroup of $G$, which is unique up to conjugacy, and let $P$ be the (normal) pro-$p$-Sylow subgroup of $I$. The quotient $I/P = T(q)$ is s split torus over the residue field, of dimension $\ell$ equal to the rank of $G$ over $k$. The Moy-Prasad filtration gives a normal subgroup $P^+$ with $P/P^+ = V$ an $F$-vector space of dimension $\ell + 1$, which decomposes into lines under the action of the torus $T$, the characters given by the $\ell + 1$ simple affine roots. 

We call a complex character $\chi$ of $P$ {\bf affine generic} if it is trivial on $P^+$ and the restriction of $\chi$ to each eigenspace for $T$ in $P/P^+$ is non-trivial. In this case, Mark showed that that the induced representation $\Ind(P|G, \chi)$ has finite length and decomposes as a direct sum of irreducible representations with multiplicity one, having distinct characters for the center $Z(q)$.

Another way to view affine generic characters is the following. The characters of $P/P^+= V$ are all of the form $\psi(\Tr(f(v))$, where $f: V \rightarrow F$ is a linear functional, $\Tr$ is the trace from $F$ to $\Z/p\Z$, and $\psi$ is a non-trivial homomorphism from $\Z/p\Z$ to $\bC^*$, and a character is affine generic if $f$ is non-zero when restricted to each eigenspace for $T$.
The ring of $T$-invariant polynomials on $V^*$ is a polynomial ring with a single generator of degree $h$, the Coxeter number of $G$. For a simple root, let $x_i$ be a basis of the $\alpha_i$ eigenspace in $V$, and let $m_i$ be the multiplicity of $\alpha_i$ in the highest root. For the negative $\alpha_0$ of the highest root, let $x_0$ be a basis of the corresponding eigenspace in $V$ and let $m_0 = 1$. Then a generating $T$-invariant polynomial of degree $h$ is given by
$$p = \prod_{i = 0}^{\ell} x_i^{m_i}.$$
The linear form $f$ on $V$ gives an affine generic character if and only if $p(f) \neq 0$. This non-vanishing condition determines the stable orbits (closed, with finite stabilizer) of $T$ on $V^*$.

We called the irreducible constituents of $\Ind(P|G, \chi)$ {\bf simple supercuspidal} representations. Mark later came up with better terminology for a larger class of compactly induced supercuspidal representations, which he discovered with Jiu-Kang Yu \cite{RY}. He called these representations {\bf epipelagic}, which refers to the first level of depth in the ocean, where sunlight can penetrate. Compared to the construction of general depth zero representations, the construction of simple supercuspidal representations requires almost nothing but the Bruhat-Tits decomposition of $G$ into double cosets for $I$. In the depth zero case, one needs the structure of all maximal compact subgroups of G and the determination of their reductive quotients, as well as the theory of Deligne and Lusztig necessary to construct discrete series representations of these reductive groups over finite fields.

\section{Simple wild parameters}

Having made such a simple construction of irreducible, supercuspidal representations, Mark and I turned to the question of their Langlands parameters. To be compatible with the conjecture of Hiraga, Ichino, and Ikeda on formal degrees, there were a number of constraints. First, the adjoint $L$-function had to be trivial, so there were no inertial invariants on the Lie algebra. This implied that the parameter was wildly ramified, so the Swan conductor of the adjoint representation was positive. The second constraint was that the Swan conductor was equal to the rank $r(G)$ -- the minimal positive value, given our condition on the invariants. We called parameters which met these two conditions {\bf simple wild parameters} \cite[\S 6]{GR}, and conjectured that they were the Langlands parameters of simple supercuspidal representations \cite[\S9]{GR}.

Simple wild parameters are fairly easy to describe, when the prime $p$ does not divide the Coxeter number $h$ of ${\bf G}$. The image of the inertia group $I$ is a finite subgroup of the normalizer of $\hat{T}$ in $\hat{G}$, isomorphic to a semi-direct product $(\Z/p\Z)^a . (\Z/h\Z)$, where $a$ is the order of $p$ modulo $h$. Wild inertia maps to a regular subgroup of $\hat{T}[p]$ and the tame quotient to a Coxeter class in the Weyl group. The corresponding ramification filtration of the inertia group $D_0$ of the finite Galois extension fixed by the kernel of the parameter is trivial at level $D_2$ \cite[\S 6.1]{GR}. A complete description of the Langlands correspondence for these packets was obtained by Kaletha \cite{Ka}.

When $p$ divides $h$, the simple wild parameters are more complicated \cite[\S 6.2]{GR}, even though the construction of simple supercuspidal representations is completely uniform. In particular, the image of inertia is not contained in the normalizer of a maximal torus $\hat{T}$.
For example, when ${\bf G} =\SL_2$ the Coxeter number $h = 2$. For primes $p > 2$ the image of inertia is a dihedral group of order $2p$ in $\PGL_2(\bC)$. However, when $p = 2$ the image of inertia is isomorphic to the alternating group $A_4$ on four letters. 

Bushnell and Henniart \cite{BH} have determined all simple wild parameters for ${\bf G} = \GL_n$, but for general groups the simple wild parameters at primes $p$ dividing $h$ remain somewhat mysterious. Reeder has discovered some exotic solvable finite subgroups in exceptional groups which are the image of these parameters. \cite[\S4]{R2}.

\section{The trace formula}

Once we had constructed simple supercuspidal representations for $p$-adic groups, a natural question to ask was whether they appeared as local components of global automorphic representations. We recall that for a reductive group over a global field $k$, with ring of ad\`eles $\A_k$, an irreducible representation  $\pi$ of  $G(\A_k)$ factors as a restricted tensor product $\pi = \otimes' \pi_v$ with almost all local components $\pi_{\frak p}$ unramified. For the purpose of this paper, we will call $\pi$ automorphic if it appears as a discrete summand of the representation of $G(\A_k)$ by right translation on a suitable space of functions on the coset space $G(k) \backslash G(\A_k)$, so has a linear form invariant under the subgroup $G(k)$.

When $G = \PGL_2$, the local component $\pi_p$ of an ad\`elic representation $\pi$ is simple supercuspidal if and only if the power of $p$ in the conductor of the standard two dimensional representation is equal to $3$. The Fermat cubic $x^3 + y^3 = 1$ is modular, and corresponds to a holomorphic modular form of weight $2$ and level $27 = 3^3$. This gives an example of an automorphic representation $\pi$ of $\PGL_2$ whose local component $\pi_3$ is simple supercuspidal. But we wanted something more systematic. For this, we turned to the trace formula.

The trace formula greatly simplifies if one works over a totally real number field, and assumes that the local components of $G$ at the real places are all compact. Here I will summarize what happens over $\Q$. Let $\A$ be the ring of ad\`eles of $\Q$ and let $G$ be a simple, simply-connected group over $\Q$ with $G(\R)$ compact. Then $G(\Q)$ is a discrete and co-compact subgroup of $G(\A)$ and smooth, compactly supported test functions on $G(\A)$ have a trace on the space of automorphic forms, which can be computed as a sum of orbital integrals over conjugacy classes in $G(\Q)$. Let $S$ denote a finite set of primes, which is non-empty and contains all of the primes $p$ where the group $G$ is not quasi-split and split by an unramified extension of $\Q_p$, and fix a (finite dimensional) irreducible representation $V$ of $G(\R)$. David Pollack and I had previously used the trace formula to count the number of automorphic representations $\pi$ with $\pi_{\infty} \cong V$, $\pi_p$ isomorphic to the Steinberg representation of $G(\Q_p)$, for primes $p \in S$, and $\pi_p$ unramified for all primes $p$ which do not lie in $S$ \cite{GP}.
Results of Bob Kottwitz on the Euler-Poincar\'e function \cite{K} allowed us to rewrite the sum of orbital integrals as a sum of stable orbital integrals over the stable torsion conjugacy classes in $G(\Q)$. The primary contribution the stable orbital integral of $\gamma$ was the term
$$\frac{1}{2^{\ell}} L_S(M_{\gamma},0) \times \Tr(\gamma|V^*).$$
Here $\ell$ is the dimension of a maximal compact torus in $G(\R)$, $M_{\gamma}$ is the Artin-Tate motive \cite{Gr2} associated to the reductive centralizer $G_{\gamma}$ (which is well-defined up to inner twisting), and 
$$L_S(M_{\gamma},s) = \prod_{p\notin S} (\det(1 - \Fr(p) p^{-s}|M_{\gamma})^{-1}$$
is its Artin $L$-function, with the primes in $S$ removed from the Euler product. The Euler product defining $L_S(M_{\gamma},s)$ converges when the real part of $s$ is large, and the special value $L_S(M_{\gamma},0)$ is obtained by analytic continuation. This special value is known to be a rational number, by results of Siegel.

I decided to see if I could modify our argument to allow the local representations $\pi_p$ to be simple supercuspidal (for a fixed affine generic character $\chi_p$, and with trivial central character) at a finite set of places $T$ which was disjoint from $S$. To do this, I simply replaced the local test function in the trace formula from the characteristic function of a hyperspecial maximal compact subgroup $G(\Z_p)$ to the function $\overline{\chi_p}$ supported on the pro-$p$-subgroup $P_p$ of an Iwahori subgroup of $G(\Z_p)$. Since this was a pro-$p$-group, the only stable torsion classes with a non-zero orbital integral had order a power of $p$. Therefore, if the set $T$ contained at least two places of different residue characteristic, the only stable torsion class with a non-zero contribution was the trivial class! To calculate the orbital integral of the trivial class, I only had to multiply the previous contribution at $\gamma = 1$ by the product of the indices $(G(\Z_p)/P_p)$ for all primes $p \in T$. This gives the special value at $s = 0$ of the modified Artin $L$-function
$$L_{S,T}(M,s) = L_S(M,s) \prod_{p \in T} \det(1 - \Fr(p) p^{1-s}|M),$$
where $M = M_1$ is the motive of $G$.

In summary, when $S$ was non-empty and $T$ contained two primes of different residue characteristic, the number of automorphic representations (counted with multiplicity) which had local component $V$ at infinity, were the Steinberg representation at primes in $S$, were simple supercuspidal (for a fixed affine generic character) at primes in $T$, and were unramified outside $S$ and $T$ is essentially given by 
$$\frac{1}{2^{\ell}} L_{S,T}(M,0) \times \dim(V).$$
This is the precise answer when $G$ is both simply-connected and has trivial center (so $G$ is the anisotropic form of $G_2$, $F_4$, or $E_8$ over $\Q$). In the general case, there are some terms involving the center $\hat{Z}$ of the dual group $\hat{G}$, for primes in $S$, and the center $Z$ of $G$, for primes in $T$ \cite[\S3]{FG}.

For example, let $G$ be the anisotropic form of the group $G_2$, which is the automorphism group of Cayley's definite octonion algebra. Then the motive of $G$ is the sum
$\Q(-1) + \Q(-5)$ of Tate motives and the above expression becomes
$$\frac{1}{4} \zeta(-1) \zeta(-5) \prod_{p\in S}(1-p)(1-p^5) \prod_{p \in T}(1-p^2)(1-p^6) \times  \dim(V) \approx \frac{1}{12096} \prod_{p \in S} p^6 \prod_{p \in T}p^8 \times \dim(V).$$
Suffice it to say that there are plenty of automorphic representations of this type whose local components at places in $T$ are simple supercuspidal!

\section{Function fields}

In June of $2008$ David Roberts came through Boston. Mark and I met with him and described our results. We emphasized that the Langlands parameters of simple supercuspidal representations had no inertial invariants on the adjoint representation. Hence the adjoint $L$-function was trivial. Moreover the Swan conductor on the adjoint representation was equal to the rank $\ell$, which was minimal given the condition on invariants.

David remarked that this reminded him of some work of Pierre Deligne and Nick Katz on Kloosterman sums. Deligne had constructed an $\ell$-adic sheaf of rank $n$ and on $\mathbb G_m$ over $F = \Z/p\Z$, which was tamely ramified at $t = 0$ with monodromy a principal nilpotent element, and wildly ramified at $t = \infty$, with Swan conductor $1$.  The traces of Frobenius $F_t$ at places $t \neq 0, \infty$ of the rational function field were given by Kloosterman sums \cite{D2} . The global monodromy had been determined by Katz \cite{Ka1}, and from his work it was an exercise to check that the Swan conductor at $t =  \infty$ of the adjoint representation was equal to the rank. Katz also asked, if this $\ell$-adic sheaf corresponded to an automorphic representation. Perhaps this automorphic representation was of the type studied above.

This convinced me to reconsider the trace formula in the function field setting. Let $G$ be a split, simple, simply-connected group over the function field $k = F(X)$, where $X$ is a smooth, projective curve of genus $g$ over the finite field $F$ of order $q$. There no Archimedian places, and less was known about the orbital integrals of Kottwitz's Euler-Poincar\'e function. But assuming that his stabilization worked equally well in characteristic $p$, I could count the automorphic representations (with multiplicities) which were Steinberg at places in the non-empty set $S$, were simple supercuspidal (for a fixed affine generic character and trivial central character) at places in the non-empty disjoint set $T$, and were unramified elsewhere. The Steinberg hypothesis at $S$ allowed one to compute this as a sum of stable orbital integrals over torsion classes which were semi-simple, hence of order prime to $p$, and the simple supercuspidal hypothesis at $T$ showed that the only non-zero contributions were given by torsion classes of $p$ power order. Therefore only the trivial class contributed to the sum!

Again, the orbital integral at the trivial class is essentially given by the special value $L_{S,T}(M,0)$ of a modified Artin $L$-function, where $M$ is the motive of $G$. Since $G$ is split, $M$ is a direct sum of the Tate motives $\Q(1-d_i)$, where $d_1 = 2,d_2, \ldots, d_{\ell} = h$ are the degrees of the generating invariant polynomials on the Lie algebra $\frak{\hat{g}}$.  The special value of $L_{S,T}(M,s)$ at $s = 0$ is the product of the modified zeta function of the curve $X$ at the negative integers $1-d_i$ :
$$L_{S,T}(M,0) = \prod_{i = 1}^{\ell} \zeta_{S,T}(1 - d_i)$$
The precise formula for the number of automorphic representations, counted with multiplicities, also involves the order of the center of $\hat{G}$ for elements in $S$ and the order of the center of $G$ for elements in $T$ \cite[7.4]{Gr4} \cite[\S7]{FG}.

Weil showed that the zeta function of $k = F(X)$ has the form
$$\zeta(s) = P(\ZZ)/(1-\ZZ)(1-q\ZZ)$$
with $\ZZ = q^{-s}$ and $P(\ZZ) = 1 - a\ZZ + \ldots + q^g\ZZ^{2g}$ a polynomial of degree $2g$ with integral coefficients. Here $q$ is the order of the finite field $F$ and $g$ is the genus of the curve $X$. Hence

$$\zeta_{S.T}(s) = P(\ZZ) \times  \prod_{v\in S}(1 - \ZZ^{\deg v})/ (1 - \ZZ) \times  \prod_{w \in T}(1 - (q\ZZ)^{\deg w})/ (1-q\ZZ)$$

is a polynomial in $\ZZ =q^{-s}$ with integral coefficients, of degree 
$$D = 2g - 2 + \deg(S) + \deg(T).$$
with leading coefficient 
$$c = \pm q^{g - 1 + \deg(T)}.$$ 
Therefore the special value $\zeta_{S,T}(1-d)$ is an integer whose absolute value is approximately $c.q^{D(d-1)},$ and just as in the number field case, the integer $L_{S,T}(M,0)$ is usually very large.

In one special case however, the degree $D$ is equal to zero and $\zeta_{S,T}(s) = 1$. Namely, when the curve $X$ has genus $g = 0$ and the sets $S$ and $T$ both consist of a single place of degree one. All curves of genus $g = 0$ over the finite field $F$ are isomorphic to the projective line $\mathbb P^1$, and the automorphism group of the projective line acts triply transitively on rational points. Hence there is no loss of generality in assuming that $S = \{0\}$ and $T = \{\infty\}$. In this case, the values of $\zeta_{S,T}$ at all negative integers are all equal to $1$, and the trace formula predicts that number of these automorphic representations, counted with multiplicity, is equal to $1$. Hence an
automorphic representation $\pi$ of $G$ over $k = F(t)$ which is Steinberg at $t=0$, simple supercuspidal (for a fixed affine generic character and trivial central character) at $t = \infty$, and unramified elsewhere is predicted to be unique!

This immediately raised the question: what are the unramified local components of the representation $\pi$ at places $t \neq 0,\infty$? Via the Satake isomorphism \cite{G5}, the unramified local representation $\pi_t$ has a Langlands parameter $c(\pi_t)$, which is a semi-simple conjugacy class in the complex dual group $\hat{G}(\bC)$. Since the global representation $\pi$ is determined completely by the local representations $\pi_0$ and $\pi_{\infty}$, it has the same field of rationality as those representations. The Steinberg representation is defined over $\Q$ and any simple supercuspidal representation is defined over $\Q(\mu_p)$. Hence the conjugacy classes $c(\pi_t)$ in $\hat{G}$ should all be rational over the cyclotomic field $\Q(\mu_p)$.

\section{An irregular connection}

Via the global Langlands correspondence, the automorphic representation $\pi$ predicted to exist by the trace formula should correspond to a system of $\ell$-adic representations 
$$\phi_{\ell}: \Gal(\overline{k}/k) \rightarrow \hat{G}(\Q_{\ell}\otimes \Q(\mu_p))$$
where $k = \Z/p\Z(t)$ is the function field of the projective line over the finite field $\Z/p\Z$, $\ell$ is a rational prime not equal to $p$, and $\hat{G}$ is the dual group. Since the global correspondence is compatible with the local correspondence, the Galois
representation $\phi_{\ell}$ should be unramified outside of $\{0,\infty\}$, should be tamely ramified at the place $t = 0$ with monodromy a principal nilpotent element $N$ in $\hat{\frak g}$ and 
should be a simple wild parameter at $t = \infty$. The image of a Frobenius class $\phi_{\ell}(Fr_t)$ for $t \neq \{0,\infty\}$ should be the $\ell$-adic realization of the semi-simple conjugacy class $c(\pi_t)$. When $G = \PGL_n$ and $\hat{G} = \SL_n$, this was presumably the $\ell$-adic sheaf on $\mathbb G_m$ constructed by Deligne, where the trace of $c(\pi_t)$ was given by an $n$-variable Kloosterman sum. What the trace formula suggested was that there should exist a similar local system on $\mathbb G_m$, for any simple group!

I lectured on this question in Paris in July of $2008$, and went to lunch afterwards with Edward Frenkel. He suggested that it might be easier to construct a complex analog of the $\ell$-adic representations $\phi_{\ell}$. Namely, there should be a connection $\nabla$ on a principal $\hat{G}$ bundle on $\mathbb G_m$ over $\bC$,  which had a regular singularity at $t = 0$ with principal nilpotent monodromy, and had an irregular singularity at $t = {\infty}$ with the smallest positive slope (defined below). Edward suggested that this connection might be an Oper, in the sense of Beilinson and Drinfeld.

With that guidance, it was easy to guess the connection \cite[\S5]{FG}. Any connection $\nabla$ on the trivial $\hat{G}$ bundle is given by a $1$-form on $\mathbb G_m$ with values in the complex Lie algebra $\hat{\frak g}$. Let $N$ be a principal nilpotent element in $\hat{\frak g}$, and let $E$ be a basis of the highest root space for the opposite Borel subgroup $\hat{B}$. Consider the connection
$$\nabla = d + N(dt/t) + E(dt).$$
This has a simple pole at $t = 0$ and a double pole at $t = \infty$. The monodromy at the regular singular point $t = 0$ is the principal nilpotent element $N$. At $\infty$, the connection $\nabla$ is irregular. Note that $E$ is a representative of the minimal nilpotent orbit for $\hat{B}$ on $\frak g$.

To compute the slope at $\infty$, we recall the definition in \cite{D1}. Let $s = 1/t$ be a uniformizing parameter there. Then the slope is the positive rational number $a/b$ provided that over the tamely ramified local extension given by $u^b = s$, $\nabla$ has a pole of order $a$, with a non-nilpotent polar part. We claim that the slope is $1/h$, where $h$ is the Coxeter number of $\hat{G}$. Indeed, taking the covering given by $u^h = s$, we may write
$$\nabla = d - hN(du/u) - hE(du/u^{h+1}).$$
Making a gauge transformation by the element $g = \rho(u)$, where $\rho$ is the co-character of $\hat {T}$ given by half the sum of the positive co-roots for $\hat{B}$, we find that this connection also has the form
$$\nabla = d - h(N+E)(du/u^2) - \rho(du/u)$$
This has a pole of order $2$, with polar part $(N+E)$. The element $N+E$ is both regular and semi-simple, by a theorem of Kostant \cite[Lemma 10]{Ko}. Hence the slope at $\infty$ is equal to  $1/h$. Moreover, the element $\exp(\rho/h)$ normalizes the maximal torus centralizing $N+E$, where it acts as a Coxeter element in the Weyl group. 

We show that the image of the wild differential Galois group at $\infty$ is a Coxeter torus -- on the character group of this torus the eigenvalues for the Coxeter element are precisely the primitive $h$ roots of unity. This is a maximal torus when $\hat{G}$ is $G_2$, $F_4$ and $E_8$, and we use this, together with the regular monodromy at $t = 0$, to show that the global differential Galois group of $\mathbb G_m$ has these three complex exceptional groups as quotients.

\section{$\ell$-adic representations}

This is where matters stood until the fall of $2009$. At a Current Developments in Mathematics conference at Harvard, B\'ao-Ch\^au Ng\^o told me that, working together with Joachim Heinloth and Zhiwei Yun, he had been able to construct both the desired automorphic representation $\pi$ of $G$ over the rational function field $k = \Z/p\Z(t)$, as well as the compatible system of $\ell$-adic representations
$$\phi_{\ell}: \Gal(\overline{k}/k) \rightarrow \hat{G}(\Q_{\ell}\otimes \Q(\mu_p))$$
associated to $\pi$ by the global Langlands correspondence! \cite{HNY}. Their construction of $\pi$ was unconditional, whereas mine depended on the unproven stabilization of the trace formula in characteristic $p$. 

For simplicity, we will assume in this section that $G$ is a simple, split and simply-connected group defined over the finite field $\Z/p\Z$, so is everywhere split over the function field $k = \Z/p\Z(t)$. Their construction of the automorphic representation $\pi$ 
begins with the construction of a newform $F$ on $G(\A)$, where $\A$ is the ring of ad\`eles of $k$. More precisely, they show that there is a unique locally constant function (up to scaling) $F: G(\A) \rightarrow \bC$ with the following invariance properties. 
\begin{itemize}
\item $F$ is left invariant under $G(k)$,
\item $F$ is right invariant under the hyperspecial subgroup $G(A_t)$, for every place $t \neq 0,\infty$, 
\item $F$ is right invariant under the Iwahori subgroup $I(0) \subset G(A_0)$,
\item $F$ transforms by the affine generic character $\chi$ under the pro-$p$-group $P(\infty) \subset G(A_{\infty})$.
\end{itemize}
The unicity of $F$ is established in \cite[\S2.1]{HNY} and relies on the combinatorial structure of the double coset space
$$G(k)\backslash G(\A)/ \prod_{t \neq 0,\infty} G(A_t) = G(F[t,t^{-1}])\backslash G(F((t))) \times G(F((t^{-1}))).$$
Since this function is unique, it is an eigenfunction for the relevant local Hecke operators. Hence the $G(\A)$ translates of $F$ generate an irreducible automorphic representation, which is their construction of $\pi$.

Their construction of the system of $\ell$-adic representations associated to $\pi$ is a tour de force. To give such a homomorphism $\phi_{\ell}$ to the dual group $\hat{G}(\Q_{\ell}\otimes\Q(\mu_p))$ is equivalent to giving, for each linear representation $V_{\ell}$ of the dual group, a local system on $\mathbb G_m$ with values in $V_{\ell} \otimes \Q(\mu_p)$, in a way compatible with direct sums and tensor products. They construct these local systems using powerful techniques in geometric Langlands theory. First, they interpret $F$ as a function on the double coset space
$$G(k) \backslash G(\A) / \prod_{x \neq 0,\infty} G(A_x) \times I(0) \times P(\infty,\chi).$$
Then they view this double coset space as the $F$-rational points on the moduli stack of $G$-bundles on $\mathbb P^1$ with some level structure at $0$ and $\infty$, and convert the function $F$ to a sheaf $\mathscr A_F$ on this moduli stack \cite[\S2.2]{HNY}. Finally, they show that $\mathscr A_F$ is a Hecke eigensheaf, and that the eigenvalues of the geometrically defined Hecke operators generate the desired local systems. In particular, the images $\phi_{\ell}(Fr_t)$ for $t \neq \{0,\infty\}$ give semi-simple conjugacy classes in $\hat{G}$ over $\Q_{\ell} \otimes \Q(\mu_p)$. Since these $\ell$-adic local systems are compatible, the conjugacy class associated to $t$ is rational over $\Q(\mu_p)$. These are the Langlands parameters $c(\pi_t)$ of the local components $\pi_t$ of the automorphic representation $\pi$.

From an analysis of the ramification at $t = 0$ and $t = \infty$ they obtain, for every $\ell \neq p$, an unramified extension of $\mathbb G_m$ in characteristic $p$ whose Galois group is an compact open subgroup of the $\ell$-adic points of $G_2$, or $F_4$, or $E_8$.

\section{$F$-isocrystals}

The $\ell$-adic Galois representations constructed by Heinloth, Ngo, and Yun were formally analogous to the complex $\hat{G}$-connection that Frenkel and I had studied. But their precise mathematical relation was only established in $2019$ by Daxin Xu and Xinwen Zhu \cite{XZ}. 

They study the irregular connection $\nabla = d + N(dt/t) + (\lambda)^hE(dt)$ on $\mathbb G_m$ over the field $K = \Q_p(\mu_p) = \Q_p(\lambda)$ with $\lambda^{p-1} = -p$, and show that it has the structure of an $F$-crystal. Namely, let $A_K$ be the ring of overconvergent power series on $\mathbb P^1$ along $\infty$ with coefficients in $K$. Elements of $A_K$ are power series in one variable with coefficients in $K$ whose radius of convergence is greater than $1$. In particular, any element of $A_K$ can be evaluated on a unit in an extension field of $K$. Xu and Zhu show \cite[Thm 1.2.4]{XZ} that there is a unique element $\phi_p(x) \in \hat{G}(A_K)$ which satisfies the differential equation
$$x \frac{d\phi_p}{dx}\phi_p^{-1} + \ad_{\phi_p}(N + \lambda^hxE) = p(N + \lambda^hx^pE).$$
Namely, $\phi_p(x)$ gives a horizontal isomorphism from the Frobenius pull-back of $\nabla$ under the map $x \rightarrow x^p$ to the original $\hat{G}$ bundle with connection on $\mathbb G_m$. We emphasize that the connection $\nabla$ and its Frobenius pull-back are both algebraic, but the horizontal isomorphism $\phi_p(x)$ is analytic.

For each element $t \neq \{0,\infty]$ in $\mathbb G_m$ over the separable closure of $\Z/p\Z$, let $X$ be the Teichmuller lifting of $t$ to a root of unity in the maximal unramified extension of $\Q_p$.  Let $V$ be a representation of $\hat{G}$ over $K$, and consider the trace of the element $\prod_{i = 0}^{\deg t - 1} \phi_p(X^{p^i})$ on $V$. There is a unique semi-simple conjugacy class $\phi_p(t)$ in $\hat{G}(K)$ with this trace, for any $V$, and they prove that this is the $p$-adic realization of the Langlands parameter $c(\pi_t)$ of the unramified representation $\pi_t$.  From this $p$-adic realization, they are able to calculate the $p$-adic valuation of these Langlands parameters, and show that they are generically ordinary.

\section{Conclusion}

The Langlands correspondence offers many insights into both number theory and representation theory, allowing discoveries in one field to suggest questions in the other. After Reeder had constructed simple supercuspidal representations a natural question was the determination of their local Langlands parameters. This led to the construction of wildly ramified extensions of $p$-adic fields, as well as the discovery of interesting finite solvable subgroups of complex Lie groups. The attempt to realize these irreducible representations as local components of automorphic representations and their parameters as local components of families of compatible $\ell$-adic Galois representations led to the discovery of some remarkable local systems on the multiplicative group over $\Z/pZ$.

\def\noopsort#1{}
\providecommand{\bysame}{\leavevmode\hbox to3em{\hrulefill}\thinspace}

\end{document}